# Optimal Planning of Integrated Heat and Electricity Systems: a Tightening McCormick Approach

Lirong Deng[a], Xuan Zhang[a], Tianshu Yang[a], Hongbin Sun[a,b,*]

[a] Shenzhen Environmental Science and New Energy Technology Engineering Laboratory, Tsinghua-Berkeley Shenzhen Institute, Tsinghua Univeristy, Shenzhen, 518055, China

[b] State Key Laboratory of Power Systems, Department of Electrical Engineering, Tsinghua University, Beijing, 100084, China

[*] Corresponding author. E-mail: shb@tsinghua.edu.cn, Tel: 010-62783086, Address: Room 3-120, West Main Building, Tsinghua University, Haidian District, Beijing, China

**ABSTRACT**

In this paper, we propose a convex planning model of integrated heat and electricity systems considering variable mass flow rates. The main challenge comes from the non-convexity of the bilinear terms in the district heating network, i.e., the product of mass flow rate and nodal temperature. To resolve this issue, we first reformulate the district heating network model through equivalent transformation and variable substitution. It shows that the reformulated model has only one set of nonconvex constraints with reduced bilinear terms and the others are linear constraints. Such a reformulation not only guarantees the optimality but fastens the solving process. To relax the remaining bilinear constraints, we apply McCormick envelopes and further propose a heuristic tightening method to constrict the bounds of the McCormick approach and get a nearby feasible solution. Case studies show that the tightening McCormick method quickly solves the heat-electricity planning problem with acceptable feasibility check and optimality.

## Nomenclature

*Functions*

| | |
|---|---|
| $f_{i,t}^{\text{HB}}$ | Hour-$t$ cost function of heating boiler $i$ |
| $f_{i,t}^{\text{CHP}}$ | Hour-$t$ cost function of CHP unit $i$ |
| $f_{i,t}^{\text{TU}}$ | Hour-$t$ cost function of non-CHP thermal unit $i$ |

*Sets*

| | |
|---|---|
| $\mathcal{T}$ | Set of indices of scheduling periods |
| $\mathcal{I}^{\text{node}}$ | Set of indices of nodes in the heating network |
| In($i$) | Set of indices of nodes flowing into node $i$ in the heating network |
| Out($i$) | Set of indices of nodes flowing out of node $i$ in the heating network |
| $\mathcal{I}^{\text{bus}}$ | Set of indices of buses in the power network |
| $\mathcal{I}^{\text{line}}$ | Set of indices of lines in the power network |
| $\mathcal{I}^{\text{HB}}$ | Set of indices of heating boilers |
| $\mathcal{I}^{\text{CHP}}$ | Set of indices of CHP units |
| $\mathcal{I}^{\text{TU}}$ | Set of indices of non-CHP thermal units |
| $\mathcal{B}_i$ | Set of indices of boundary pairs in CHP unit $i$ |

*Variables*

| | |
|---|---|
| $H_i^{\text{HB}}$ | Heat output of heating boiler $i$ |
| $P_i^{\text{CHP}}$ | Active power output of CHP unit $i$ |
| $H_i^{\text{CHP}}$ | Heat output of CHP unit $i$ |

| | | | |
|---|---|---|---|
| $P_i^{TU}$ | Active power output of non-CHP thermal unit $i$ | $c$ | Specific heat capacity of water |
| $P_{Gi}$ | Active power generation at bus $i$ | $A_{i,b}^{CHP}, B_{i,b}^{CHP},$ $D_{i,b}^{CHP}$ | Parameters of boundary $b$ of the feasible region at CHP unit $i$ |
| $P_{Di}$ | Active power demand at bus $i$ | | |
| $H_{Gi}$ | Heat production at node $i$ | $\overline{P_{ij}}$ | Upper bound of the power flow across line $ij$ |
| $H_{Di}$ | Heat demand at node $i$ | | |
| $P_{ij}$ | Active power flow across line $ij$ | $L_{ji}$ | Length of pipeline from node $j$ to node $i$ |
| $\theta_{ij}$ | Phase angle difference between bus $i$ and bus $j$ | $\tau_{o,i}$ | Predetermined temperature for the source node $i$ |
| $m_{ji}$ | Mass flow rate of water transferred from node $j$ to node $i$ in the pipeline of the heating network | $\tau_a$ | Ambient temperature |
| | | $C_i^{HB}$ | Production cost coefficient of heating boiler $i$ |
| $\tau_i$ | Outlet temperature of node $i$ | $C_i^{CHP,0\sim5}$ | Generation cost coefficients of CHP unit $i$ |
| $\tau_{ji}$ | Outlet temperature of pipeline from node $j$ to node $i$ | $C_i^{TU,1\sim2}$ | Generation cost coefficients of non-CHP thermal unit $i$ |
| Parameters | | | |
| $\upsilon$ | Heat transfer coefficient per unit length | | |

## 1. Introduction

The electric power system (EPS) is a crucial part of the national energy supply, but policymakers are gradually realizing that heat supply also plays a key role in the energy system. The International Energy Agency found that over half of global energy use is for heating [1]. For heat supply, the district heating system (DHS) occupies a large proportion in high-density areas, compared with the situation that other heat alternatives, such as individual heat pumps, gas boilers, solar thermal, and electrical heating, that are often applied in low-density areas. Although different countries and areas have different proportions among those heat supply methods, the district heating approach has been proven to produce high energy efficiency [2].

Electricity and heating can be produced simultaneously with centralized energy generation and district heating infrastructures. Generally, those two large energy systems—EPS and DHS—are tightly connected by combined heat and power (CHP, also known as co-generation) plants and power-to-heat facilities. By 2050, CHP will provide 26% of electricity in the EU. In Denmark, the government aims to achieve 100% renewable heat and electricity production by 2035 [3]. It is expected that the EPS and the DHS will impact each other to a greater extent in the energy production and consumption process in the near future. Therefore, the idea of integrated heat and electricity system has attracted interests from practitioners and researchers.

Studies on integrated heat and electricity systems are emerging in different aspects, including modeling [4], state estimation [5], unit commitment [6], economic dispatch [7, 8, 9], market mechanisms [10, 11], planning [12, 13], etc. Among them, modeling plays a fundamental and substantial role. Although EPS modeling has been studied thoroughly during the past decades, there are still many problems in the modeling of DHS. Generally, DHS has three regulation modes, shown in Table I. Quality regulation employs a constant mass flow rate and variable temperature strategy [15]. Hence, the constraints in related optimization problems become linear and easy to handle. In contrast, quantity regulation maintains a constant supply temperature but regulates the quantity of mass flow rates. It brings flexibility and cost reduction because mass flow rates will change according to the heat load, which will reduce the power consumption of circulating water pumps. Apparently, more economic efficiency and flexibility can be achieved by regulating both temperatures and mass flow rates, which is referred to as quality-quantity regulation. We argue that from the computation perspective, complexities of quantity and quality-quantity regulations are almost the same. The reason is that in the quantity regulation, except for the source nodes, nodal temperatures are indeed variable. Temperature will drop along the transfer process because of heat loss, where heat loss is further affected by variable mass flow rates. In this paper, we will focus on quantity regulation, but the method can be directly applied to the quality-quantity regulation.

Table I. Different regulation modes in DHS

|  | Mass flow rate | Supply nodal temperature |
|---|---|---|
| Quality regulation | constant | variable |
| Quantity regulation | variable | constant |
| Quality-quantity regulation | variable | variable |

Modeling of quantity-regulated DHS belongs to the pooling problem [14]. The pooling problem is a network flow problem that aims to find the minimum cost way to mix several inputs in intermediate pools such that the output meets the demand or certain requirements. To mix inputs is to mix the product of the flow quantity and feature. As a result, the pooling problem becomes a bilinear program. In the DHS modeling, bilinear terms come from the product of mass flow rate and temperature. Current achievements to deal with bi-linearity of the DHS optimization include:

- Gradient-based methods, such as interior point method conducted in IPOPT. Those methods are easy to implement, but present local minimum and slow convergence when the network becomes large.
- Generalized Benders Decomposition [15, 16]. It can solve certain types of nonlinear programming and mixed-integer nonlinear programming, but sometimes cannot guarantee to converge to the global optimum.
- Relaxations, such as conic relaxation [17] and polyhedral relaxation [18]. Those methods are convex thus can be solved efficiently by off-the-shelf solvers, but may lead to infeasible solutions to the original problem.
- Relaxation tightening methods, see [14, 18, 19] for example. Among them, branch and bound, or branch and cut, has been successfully implemented in the latest version of Gurobi, i.e., Gurobi 9.0 [20], to deal with bilinear programming. These methods guarantee global optimum and can be used as a benchmark to evaluate the optimality of other methods. However, when dealing with large-scale problems, it may converge slowly.

In this paper, our goal is to propose a convex model and an efficient algorithm for the integrated heat and electricity systems, where the algorithm is expected to find the global optimum or a near-global-optimal feasible solution with relatively small computational burden. To begin with, we reformulate the classical quantity-regulated DHS optimization model, through equivalent transformation and first-order Tylor expansion. Compared with the original model which has two bilinear terms in each non-convex constraint, the reformulated model reduces the bilinear terms by half. Such a reformulation not only guarantees the optimality but reduces the computational complexities of the original problem. Then, we perform McCormick envelopes to convexify the bilinear constraints and get a lower bound of the reformulated model. Since the lower bound may not be feasible, a heuristic tightening method is further established to constrict the bounds of the McCormick approach and get a nearby feasible solution. Notably, in this paper, we apply the proposed model and algorithm into the energy production planning model with daily or weekly time scales. However, the proposed methods can also be extended to intra-day or real-time economic dispatch of integrated systems, if the newly added thermal dynamics can be formulated as convex or bilinear constraints. Bilinear-form of thermal dynamics have been realized in [18]. Future work could consider the thermal dynamics and conduct the economic dispatch analysis with the proposed methods.

In summary, our contributions are mainly three-fold: 1) the classic integrated heat and electricity planning problem with quantity-regulated DHS model is reformulated through variable substitution and equivalent transformation, which largely reduces the bilinear complexity of the classic model. 2) Based on the reformulated model, McCormick envelopes are applied to relax the remaining bilinear terms. Furthermore, a tightening method is proposed to constrict the upper and lower bounds of the McCormick envelopes through perturbation of the latest optimal results. 3) Case studies show that the tightening McCormick method quickly solves the problem with acceptable feasibility check and optimality. Meanwhile, with the convex property, the tightening McCormick method is promising in large-scale integrated heat and electricity optimization and can allow economic analysis using shadow prices.

The remainder of this paper is organized as follows. Section II introduces the integrated planning problem with a DHS based model and a reformulated model. Section III details the convex relaxation and the tightening McCormick algorithm to solve the integrated problem. Section IV performs case studies and section VII draws conclusions.

## 2. Problem Formulation

The planning model consists of the DHS and the EPS. For DHS, we present a *base* model first, which is a nonlinear optimization problem without any relaxation of the constraints. Then, we derive an equivalent *reformulated* model of DHS through first-order Taylor expansion and variable substitution. The reformulated model turns out to reduce the nonlinearity of the base model. For EPS, we adopt the direct current (DC) power flow model. Coupling elements connecting the two systems are CHP units. The planning model determines the heat and electricity production of the generation units so as to minimize the total cost and meet operation constraints. For notational simplicity, we build single-time horizon models for DHS and EPS, respectively. Then, we extend the whole planning optimization to the multi-time horizon case and add the notation of time index $t$ to all the decision variables.

### 2.1 District Heating System Modeling

In the quantity-regulated radial district heating network, we regulate mass flow rate through the circulating pump but fix the supply temperature of heat sources. Mass flow rate has magnitude and direction. In this paper, we assume the magnitude of mass flow rate is variable, but the direction is fixed. It is reasonable because frequent change of direction will lead to supply instability.

**Base Model**
$$H_{Gi} - H_{Li} + c \sum_{j \in \text{In}(i)} m_{ji} \tau_{ji} = c \sum_{k \in \text{Out}(i)} m_{ik} \tau_i, \forall i \in \mathcal{I}^{\text{node}}, \tag{1a}$$

$$\sum_{j \in \text{In}(i)} m_{ji} = \sum_{k \in \text{Out}(i)} m_{ik}, \forall i \in \mathcal{I}^{\text{node}}, \tag{1b}$$

$$\tau_{ji} = \left(\tau_j - \tau_a\right)\exp\left(-\frac{\upsilon L_{ji}}{cm_{ji}}\right) + \tau_a, \forall ji \in \mathcal{I}^{\text{pipe}}, \tag{1c}$$

$$\tau_i^{\min} \leq \tau_i \leq \tau_i^{\max}, \tau_{ji}^{\min} \leq \tau_{ji} \leq \tau_{ji}^{\max}, \forall i \in \mathcal{I}^{\text{node}} \setminus \left\{\mathcal{I}^{\text{CHP}} \bigcup \mathcal{I}^{\text{EB}}\right\}, \tag{1d}$$

$$m_{ji}^{\min} \leq m_{ji} \leq m_{ji}^{\max}, \forall ji \in \mathcal{I}^{\text{pipe}}, \tag{1e}$$

$$\tau_i = \tau_{o,i}, \forall i \in \mathcal{I}^{\text{CHP}} \bigcup \mathcal{I}^{\text{EB}}, \tag{1f}$$

where constraints (1a) define the nodal heat balance. Constraints (1b) are nodal flow balance. Constraints (1c) describe the transferring process of water, considering heat-loss factors. Constraints (1d) are minimum- and maximum-operating limits of the nodal outlet temperature and pipeline outlet temperature. Constraints (1e) gives the minimum- and maximum-operating limits of mass flow rates. Constraints (1f) are the predetermined temperatures for source nodes. Note that eliminating (1f) from (1) constitutes the quality-quantity regulation model, whose computational complexity keeps the same as the quantity-regulated model.

### 2.2 DHS Reformulation

Non-convexity of the DHS base model (1) comes from (1a) and (1c). One is the bilinear term $m\tau$. The other is the exponential function $\exp(-\upsilon L/(cm))$.

- To deal with the bilinear term, we introduce auxiliary variable $H=m\tau$. It turns out that with this variable substitution, bilinear terms are reduced.
- To deal with the exponential function term, we use the first order Taylor expansion to approximate constraint (1c).

Denote $\tilde{\tau}_j = \tau_j - \tau_a$, $\tilde{\tau}_{ji} = \tau_{ji} - \tau_a$, thus (1a) and (1c) are equivalent to

$$H_{Gi} - H_{Li} + c \sum_{j \in \text{In}(i)} \left(m_{ji}\tilde{\tau}_{ji} + m_{ji}\tau_a\right) = c \sum_{k \in \text{Out}(i)} \left(m_{ik}\tilde{\tau}_i + m_{ik}\tau_a\right), \quad (2)$$

$$\tilde{\tau}_{ji} = \tilde{\tau}_j \cdot \exp\left(-\frac{\upsilon L_{ji}}{cm_{ji}}\right). \quad (3)$$

By introducing auxiliary variables

$$H_{ji}^{\text{in}} = cm_{ji}\tilde{\tau}_{ji}, \quad (4)$$

$$H_{ik}^{\text{out}} = cm_{ik}\tilde{\tau}_i, \quad (5)$$

we transform (2) and (3) into

$$H_{Gi} - H_{Li} + \sum_{j \in \text{In}(i)} H_{ji}^{\text{in}} = \sum_{k \in \text{Out}(i)} H_{ik}^{\text{out}}, \forall i \in \mathcal{I}^{\text{node}}, \quad (6)$$

$$H_{ji}^{\text{in}} = H_{ji}^{\text{out}} \cdot \exp\left(-\frac{\upsilon L_{ji}}{cm_{ji}}\right). \quad (7)$$

Since normally $\upsilon L_{ji} \ll cm_{ji}$, we can use the first order Taylor expansion to approximate (7).

$$H_{ji}^{\text{in}} = H_{ji}^{\text{out}} \cdot \exp\left(-\frac{\upsilon L_{ji}}{cm_{ji}}\right) \approx H_{ji}^{\text{out}} \cdot \left(1 - \frac{\upsilon L_{ji}}{cm_{ji}}\right) = H_{ji}^{\text{out}} - \upsilon L_{ji}\tilde{\tau}_j, \forall ji \in \mathcal{I}^{\text{pipe}}. \quad (8)$$

To transform upper and lower limits constraints associated with temperature (i.e., (1d)) into constraints related to heat power $H$, we have the following:

$$cm_{ji}\left(\tau_{ji}^{\min} - \tau_a\right) \leq H_{ji}^{\text{in}} \leq cm_{ji}\left(\tau_{ji}^{\max} - \tau_a\right), \forall ji \in \mathcal{I}^{\text{pipe}}, \quad (9)$$

$$cm_{ik}\left(\tau_i^{\min} - \tau_a\right) \leq H_{ik}^{\text{out}} \leq cm_{ik}\left(\tau_i^{\max} - \tau_a\right), \forall ik \in \mathcal{I}^{\text{pipe}}. \quad (10)$$

Therefore, the district heating network formulation under quantity regulation (1) can be reformulated as (11), which is termed as reformulated model.

**Reformulated Model**
$$H_{Gi} - H_{Li} + \sum_{j \in \text{In}(i)} H_{ji}^{\text{in}} = \sum_{k \in \text{Out}(i)} H_{ik}^{\text{out}}, \forall i \in \mathcal{I}^{\text{node}}, \quad (11a)$$

$$H_{ji}^{\text{in}} = H_{ji}^{\text{out}} - \upsilon L_{ji}\tilde{\tau}_j, \forall ji \in \mathcal{I}^{\text{pipe}}, \quad (11b)$$

$$cm_{ji}\left(\tau_{ji}^{\min} - \tau_a\right) \leq H_{ji}^{\text{in}} \leq cm_{ji}\left(\tau_{ji}^{\max} - \tau_a\right), \forall ji \in \mathcal{I}^{\text{pipe}}, \quad (11c)$$

$$cm_{ik}\left(\tau_i^{\min} - \tau_a\right) \leq H_{ik}^{\text{out}} \leq cm_{ik}\left(\tau_i^{\max} - \tau_a\right), \forall ik \in \mathcal{I}^{\text{pipe}}, \quad (11d)$$

$$\tau_i^{\min} - \tau_a \leq \tilde{\tau}_i \leq \tau_i^{\max} - \tau_a, \forall i \in \mathcal{I}^{\text{node}} \setminus \{\mathcal{I}^{\text{CHP}} \cup \mathcal{I}^{\text{HB}}\}, \quad (11e)$$

$$H_{ik}^{\text{out}} = cm_{ik}\tilde{\tau}_i, \quad (11f)$$

$$(1b), (1e), (1f). \quad (11g)$$

Note that (4) does not appear in (11). It has been eliminated because $H_{ji}^{\text{in}}$ totally represents $\tilde{\tau}_{ji}$. Introducing $H_{ji}^{\text{in}}$ and eliminating variables $\tilde{\tau}_{ji}$ and (4) will not influence the feasible region of the original problem. However, (11f) is not with the same case. The reason is that $\tilde{\tau}_i$ is not only in (11f), but also appears in (11b). In the reformulated model (11), all constraints are linear except for (11f). Arguably, (11) is equivalent to (1), with negligible errors.

**2.3 Electric Power System Modeling**

$$P_{Gi} - P_{Li} = \sum_{j \in \mathcal{I}_i^{\text{bus}}, j \neq i} \frac{\theta_i - \theta_j}{X_{ij}}, \forall i \in \mathcal{I}^{\text{bus}}, \quad (12a)$$

$$-\bar{P}_{ij} \leq \frac{\theta_i - \theta_j}{X_{ij}} \leq \bar{P}_{ij}, \forall ij \in \mathcal{I}^{\text{line}}. \tag{12b}$$

Constraints (12a) provide the nodal active power balance equations. Constraints (12b) enforce limits of transmission lines.

**2.4 Energy Sources Modeling**

There are three typical types of energy sources in the integrated market: heating boilers, CHP units, and non-CHP thermal units. Usually, the generation cost of a CHP unit is formulated as a convex quadratic cost function. As for heating boilers, the generation cost is a linear function. The objectives of the above energy sources are described as follows

$$f_i^{\text{HB}} = C_i^{\text{HB}} \cdot H_i^{\text{HB}}, \forall i \in \mathcal{I}^{\text{HB}}, \tag{13a}$$

$$\begin{aligned} f_i^{\text{CHP}} &= C_i^{\text{CHP},0} + C_i^{\text{CHP},1} \cdot P_i^{\text{CHP}} + C_i^{\text{CHP},2} \cdot \left(P_i^{\text{CHP}}\right)^2 \\ &+ C_i^{\text{CHP},3} \cdot H_i^{\text{CHP}} + C_i^{\text{CHP},4} \cdot \left(H_i^{\text{CHP}}\right)^2 + C_i^{\text{CHP},5} \cdot P_i^{\text{CHP}} \cdot H_i^{\text{CHP}}, \forall i \in \mathcal{I}^{\text{CHP}}, \end{aligned} \tag{13b}$$

$$f_i^{\text{TU}} = C_i^{\text{TU},2} \cdot \left(P_i^{\text{TU}}\right)^2 + C_i^{\text{TU},1} \cdot P_i^{\text{TU}}, \forall i \in \mathcal{I}^{\text{TU}}. \tag{13c}$$

The following constraints impose operational bounds for those energy sources, respectively. For CHP units, bounds usually refer to the boundaries of the feasible region, whose shape may be either line or polygon [10], representing the relations between heat output and electricity output as well as their upper/lower limits.

$$\underline{H_i^{\text{HB}}} \leq H_i^{\text{HB}} \leq \overline{H_i^{\text{HB}}}, \forall i \in \mathcal{I}^{\text{HB}}, \tag{14a}$$

$$A_{i,b}^{\text{CHP}} \cdot P_i^{\text{CHP}} + B_{i,b}^{\text{CHP}} \cdot H_i^{\text{CHP}} \leq D_{i,b}^{\text{CHP}}, \forall i \in \mathcal{I}^{\text{CHP}}, b \in \mathcal{B}_i, \tag{14b}$$

$$\underline{P_i^{\text{TU}}} \leq P_i^{\text{TU}} \leq \overline{P_i^{\text{TU}}}, \forall i \in \mathcal{I}^{\text{TU}}. \tag{14c}$$

**2.5 Planning Problem**

The objective function minimizes the total generation costs of all generation units over the *T*-hour model horizon. Constraints include (1) for base model or (11) for reformulated model, (12), (14) and two other constraints regarding nodal electricity/heat production equalities:

$$\min_x \; F = \sum_{t \in \mathcal{T}} \left( \sum_{j \in \mathcal{I}^{\text{TU}}} f_{j,t}^{\text{TU}} + \sum_{j \in \mathcal{I}^{\text{CHP}}} f_{j,t}^{\text{CHP}} + \sum_{j \in \mathcal{I}^{\text{HB}}} f_{j,t}^{\text{HB}} \right)$$

$$x := \left\{ P_{j,t}^{\text{TU}}, P_{j,t}^{\text{CHP}}, H_{j,t}^{\text{CHP}}, H_{j,t}^{\text{HB}}, V_{i,t}, \theta_{i,t}, \tau_{i,t}, \tau_{ji,t} \right\} \tag{15a}$$

$$\text{s.t.} \quad (1) \text{ or } (11), (12), (14), \forall t \in \mathcal{T}, \tag{15b}$$

$$P_{Gi,t} = \sum_{j \in \mathcal{I}_i^{\text{CHP}}} P_{j,t}^{\text{CHP}} + \sum_{j \in \mathcal{I}_i^{\text{TU}}} P_{j,t}^{\text{TU}}, \forall i \in \mathcal{I}^{\text{bus}}, t \in \mathcal{T}, \tag{15c}$$

$$H_{Gi,t} = \sum_{j \in \mathcal{I}_i^{\text{HB}}} H_{j,t}^{\text{HB}} + \sum_{j \in \mathcal{I}_i^{\text{CHP}}} H_{j,t}^{\text{CHP}}; \forall i \in \mathcal{I}^{\text{node}}, t \in \mathcal{T}, \tag{15d}$$

where $\mathcal{I}_i^{\text{HB}}$, $\mathcal{I}_i^{\text{CHP}}$, $\mathcal{I}_i^{\text{TU}}$ are sets of heating boilers, CHP units, non-CHP thermal units that are connected to bus *i*, respectively. The problem (15), whether with DHS base model (1) or reformulated model (11), is a bilinear program. It is NP-hard.

## 3. Convex Relaxation and Solution Algorithm of Integrated Heat and Electricity Problem

In the DHS reformulated model, the bilinear term is the only difficulty. An intuitive idea is to remove it, i.e., (11f), directly to check whether the bilinear constraints have a great impact on the quality of the solution. However, simulation shows removing bilinear terms leads to solutions with large violation errors to the constraints (11f). A classic method to relax bilinear terms is McCormick envelopes. The optimal results from the relaxed McCormick model are global optimum but may not be feasible to the reformulated model. Hence, a tightening McCormick algorithm has been proposed to tighten the boundaries of the McCormick relaxation and find a feasible result near the global optimum.

### 3.1 McCormick Convex Relaxation

We replace (11f) with the McCormick relaxation.

**McCormick Model** (11a)-(11e), (11g) and (16a)

$$H_{ik}^{out} \geq c\left(m_{ik}^{min}\tilde{\tau}_i + \tilde{\tau}_i^{min}m_{ik} - m_{ik}^{min}\tilde{\tau}_i^{min}\right) \tag{16b}$$

$$H_{ik}^{out} \geq c\left(m_{ik}^{max}\tilde{\tau}_i + \tilde{\tau}_i^{max}m_{ik} - m_{ik}^{max}\tilde{\tau}_i^{max}\right) \tag{16c}$$

$$H_{ik}^{out} \leq c\left(m_{ik}^{max}\tilde{\tau}_i + \tilde{\tau}_i^{min}m_{ik} - m_{ik}^{max}\tilde{\tau}_i^{min}\right) \tag{16d}$$

$$H_{ik}^{out} \leq c\left(m_{ik}^{min}\tilde{\tau}_i + \tilde{\tau}_i^{max}m_{ik} - m_{ik}^{min}\tilde{\tau}_i^{max}\right) \tag{16e}$$

By convexifying the problem, a global minimum is achieved in the planning problem. This global minimum can be regarded as a lower bound to the original problem. However, the McCormick relaxation still induces relatively large errors of the bilinear constraints.

### 3.2 Tightening McCormick Algorithm

Since we have found an optimal solution through the McCormick model, then to reduce errors of the bilinear constraints, a feasible solution near the optimal solution is expected. A basic idea is that if we can get tighter bounds of the McCormick method, we may approach to the near-optimal results that less violate the bilinear constraints. As such, we employ a simple bound-tightening algorithm to improve the McCormick method, shown in Algorithm 1. The main idea of Algorithm 1 is that at each iteration $n$, the upper and lower bounds of decision variables, i.e., mass flow rates and nodal temperature, are updated according to the results from the last iteration $n-1$ and a sequence of hyperparameter $\varepsilon$, $0<\varepsilon<1$. The principle to choose the sequence is to gradually reduce the value of $\varepsilon$ so as to tighten the bounds [21]. Moreover, to ensure feasibility of the original problem, the updated bounds should be the intersection of the result-orientated bounds and the original bounds.

---

**Algorithm 1:** Tightening McCormick

1: **inputs:** $\delta$, $N$, $\{\varepsilon^n\}_{n=1}^{N}$

2: $n \leftarrow 1$; $m_{ik}^{min,ini} \leftarrow m_{ik}^{min}$, $m_{ik}^{max,ini} \leftarrow m_{ik}^{max}$, $\forall ik \in \mathcal{I}^{pipe}$; $\tilde{\tau}_i^{min,ini} \leftarrow \tilde{\tau}_i^{min}$, $\tilde{\tau}_i^{max,ini} \leftarrow \tilde{\tau}_i^{max}$, $\forall i \in \mathcal{I}^{node} \setminus \{\mathcal{I}^{CHP} \cup \mathcal{I}^{HB}\}$

3: **while** $n \leq N$ **do**

4:     Solve (16) and get $H_{ik}^{out}$, $m_{ik}$, $\tilde{\tau}_i$

5:     **if** $\left|H_{ik}^{out} - cm_{ik}\tilde{\tau}_i\right|/H_{ik}^{out} \leq \delta$ **then**

6:        **break;**

7:     **end**

8:     $m_{ik}^{min} \leftarrow \max\left\{(1-\varepsilon^n)m_{ik}, m_{ik}^{min,ini}\right\}$, $\forall ik \in \mathcal{I}^{pipe}$

9:     $m_{ik}^{max} \leftarrow \min\left\{(1+\varepsilon^n)m_{ik}, m_{ik}^{max,ini}\right\}$, $\forall ik \in \mathcal{I}^{pipe}$

10:    $\tilde{\tau}_i^{min} \leftarrow \max\left\{(1-\varepsilon^n)\tilde{\tau}_i, \tilde{\tau}_i^{min,ini}\right\}$, $\forall i \in \mathcal{I}^{node} \setminus \{\mathcal{I}^{CHP} \cup \mathcal{I}^{HB}\}$

11:     $\tilde{\tau}_i^{\max} \leftarrow \min\left\{\left(1+\varepsilon^n\right)\tilde{\tau}_i, \tilde{\tau}_i^{\max,\text{ini}}\right\}$, $\forall i \in \mathcal{I}^{\text{node}} \setminus \left\{\mathcal{I}^{\text{CHP}} \cup \mathcal{I}^{\text{HB}}\right\}$
12:     $n \leftarrow n+1$
13: **end while**

## 4. Case Studies

In this section, we compare the performance of several models, shown in Table II:

**Table II.** DHS Models settings for comparison

| Model | Formulation | Feature | Solver | Optimum |
|---|---|---|---|---|
| Base | (1) | Nonconvex | Bilinear solver in Gurobi | Global |
|  |  |  | IPOPT | Local |
| Reformulated | (11) | Nonconvex | Bilinear solver in Gurobi | Global |
| Remove-bilinear | (11) without (11f) | Convex | Convex solver in Gurobi | Global |
| McCormick | (16) | Convex | Convex solver in Gurobi | Global |
| Tightening McCormick | Algorithm 1 | Convex | Convex solver in Gurobi | Global |

We test those models in two integrated heat and electricity systems: one is a small-scale system with a 6-bus EPS and an 8-node DHS (shown in Fig. 1); the other is a large-scale system with a 118-bus EPS and a 33-node DHS. Detailed data can be found in [22].

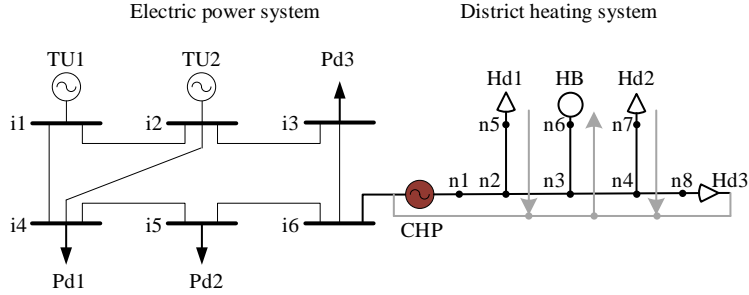

**Fig. 1.** An integrated heat and electricity system with a 6-bus EPS and an 8-node DHS.

### 4.1. Optimality

Table III shows the objective-function values and gaps in two cases, respectively. Base and reformulated models solved by the bilinear solver in Gurobi 9.0 provide the global optimal results. Note that it shows there is almost no gap between those two models, proving that the reformulation preserves the same solution structure. As expected, the remove bilinear model gives the lowest bound of the optimal solution while McCormick model provides a tighter lower bound. Base model solved by the local solver IPOPT presents an upper bound of the optimal solution. Tightening McCormick model has a small value gap and fast solution time. When the network becomes large, the base model with the bilinear solver in Gurobi is time-consuming to solve. Since the problem is highly non-convex, IPOPT even fails to converge in the large-scale case. However, with the convexified model, the tightening McCormick algorithm can be applied to the large-scale integrated system optimization with relatively small errors and efficiently solved by the off-the-shelf solvers. In all, the tightening McCormick model performs well on both solution accuracy and computational efficiency. Moreover, compared to using global search like branch and bound in Gurobi 9.0, the tightening McCormick model is a convex model and it is easy to derive shadow prices for further economic analysis.

From the perspective of the regulation methods, we can observe that regulation with variable mass flow rate has lower costs than that with constant mass flow rate in Table III, since variable mass flow rate case results in more flexibility in both mass flowrate rate and temperature. Detailed day-ahead scheduling has been demonstrated in Fig. 2.

Table III. Objective-function values and gaps

(a) the small-scale case

| | Model | Value | Gap[1] | Solver time(s) |
|---|---|---|---|---|
| | Base (Global) | 125381.24 | / | 0.10 |
| | Base (Local) | 125383.75 | 0.0020% | 196.30 |
| Variable mass flow rate | Reformulated | 125381.69 | 0.0004% | 1.83 |
| | Remove Bilinear | 125250.57 | 0.1042% | 0.03 |
| | McCormick | 125356.71 | 0.0196% | 0.12 |
| | Tightening McCormick | 125364.34 | 0.0135% | 0.42 |
| Constant mass flow rate | / | 125839.55 | / | / |

(b) the large-scale case

| Regulation | Model | Value | Gap[1] | Solver time(s) |
|---|---|---|---|---|
| | Base (Global) | 118005.81 | / | 9472.53 |
| | Base (Local) | -[2] | -[2] | >9999 |
| Variable mass flow rate | Reformulated | 118005.95 | 0.0001% | 4463.66 |
| | Remove Bilinear | 117997.82 | 0.0069% | 0.59 |
| | McCormick | 118001.31 | 0.0039% | 1.14 |
| | Tightening McCormick | 118005.47 | 0.0004% | 7.42 |
| Constant mass flow rate | | / | 118049.62 | / | / |

1 Gap is defined as the relative error with respect to the value from Base (Global).
2 - no solution found within 9999s.

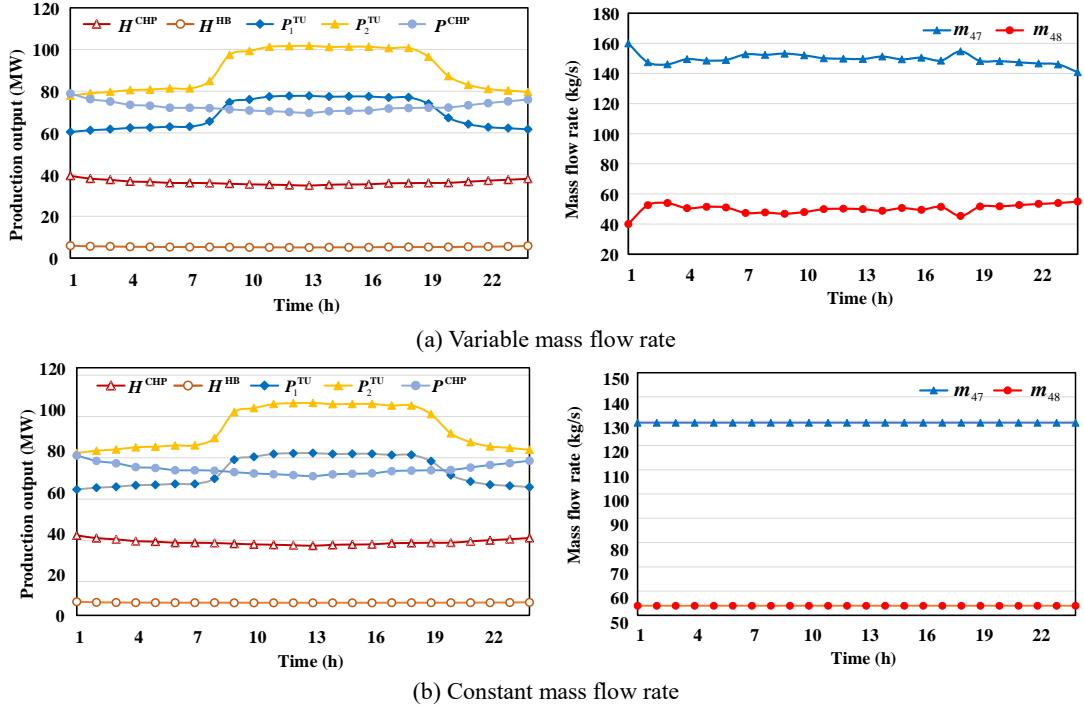

(a) Variable mass flow rate

(b) Constant mass flow rate

**Fig. 2.** Day-ahead scheduling of the small-scale case with different regulation methods.

### 4.2. Feasibility

Fig. 3 shows violations of constraints (11f) in two cases, respectively. Violations are represented by the errors of the bilinear constraints, i.e., $\left|H_{ik}^{out} - cm_{ik}\tilde{\tau}_i\right|/H_{ik}^{out}$. We can observe that although the remove bilinear and McCormick models provide lower costs, their violations are huge, which cannot be acceptable when applied to the real world operations. However, the tightening McCormick model gives relatively small violation errors. Maximum errors are no more than 7.819% and average errors are less than 2.231% in all cases. Detailed error data can be found in Table IV.

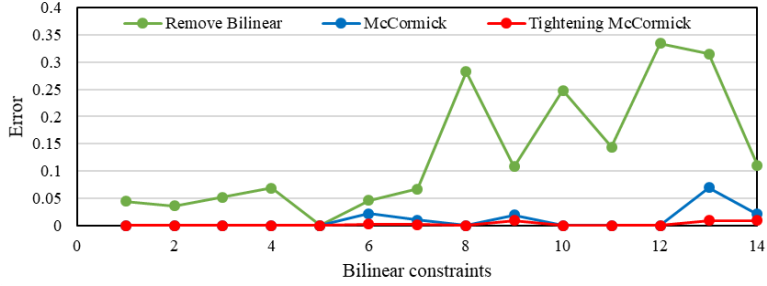

(a) the small-scale case

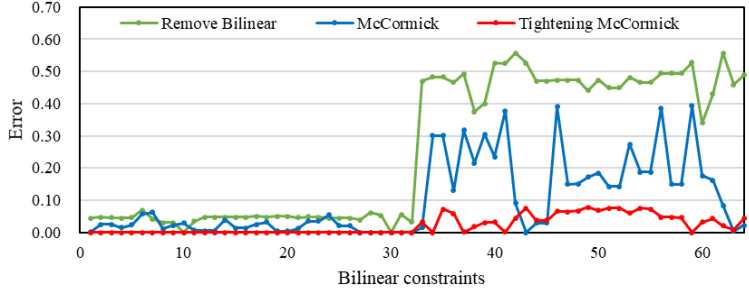

(b) the large-scale case

**Fig. 3.** Violations of constraints (11f) at hour 1

**Table IV**. Errors of the bilinear constraints (11f)

(a) the small-scale case

|  | Remove Bilinear | McCormick | Tightening McCormick |
| --- | --- | --- | --- |
| Average Errors | 13.269% | 1.027% | 0.239% |
| Maximum Errors | 33.445% | 7.000% | 0.967% |

(b) the large-scale case

|  | Remove Bilinear | McCormick | Tightening McCormick |
| --- | --- | --- | --- |
| Average Errors | 25.869% | 10.084% | 2.231% |
| Maximum Errors | 55.679% | 39.303% | 7.819% |

## 6. Conclusions

In this paper, we have proposed a convex planning model of the integrated heat and electricity systems. To reduce the bilinear terms coming from the quantity-regulated heat flow, we first reformulate the model through variable substation and equivalent transformation. Then, the remaining bilinear terms are relaxed by the McCormack envelopes. To guarantee the feasibility of the results, a bound-tightening method is presented to improve the bounds of the McCormick method and get a nearby feasible solution. Case studies show that the tightening McCormick method quickly solves the problem with acceptable feasibility check and optimality. Meanwhile, with convexity, the tightening McCormick method is promising in large-scale integrated heat and electricity system optimization and can allow economic analysis using shadow prices. Note that the performance of tightening McCormick algorithm highly depends on the chosen sequence $\{\varepsilon^k\}_{k=1}^{K}$. Further studies could focus on finding efficient sequence and perform parameter analysis.

## Acknowledgements

This work was supported by Science and Technology Program of State Grid Corporation of China (522300190008).

## References

1. IEA. Renewable Heat Policies, IEA, Paris, 2018, https://www.iea.org/reports/renewable-heat-policies.
2. State of Green. District Energy. Denmark, 2016.
3. PA Østergaard, H Lund. A renewable energy system in Frederikshavn using low-temperature geothermal energy for district heating. Appl. Energy, 2011, 88(2): 479-487.
4. X. Liu, J. Wu, N. Jenkins, A. Bagdanavicius. Combined analysis of electricity and heat networks. Appl. Energy, 2016, 162: 1238-1250.
5. T. Sheng, Q. Guo, H. Sun, Z Pan, J Zhang. Two-stage state estimation approach for combined heat and electric networks considering the dynamic property of pipelines. Energy Procedia, 2017, 142: 3014-3019.
6. Z. Li, W. Wu, J. Wang, B. Zhang, T. Zheng. Transmission-constrained unit commitment considering combined electricity and district heating networks. IEEE Trans. on Sustain. Ener., 2015, 7(2): 480-492.
7. C. Lin, W. Wu, B. Zhang, and Y. Sun. Decentralized solution for combined heat and power dispatch through benders decomposition. IEEE Trans. on Sustain. Ener., 2017, 8(4): 1361-1372.
8. Z. Pan, Q. Guo, and H. Sun. Feasible region method based integrated heat and electricity dispatch considering building thermal inertia. Appl. Energy, 2017, 192: 395-407.
9. Y. Xue, Z. Li, C. Lin, Q. Guo, H. Sun. Coordinated dispatch of integrated electric and district heating systems using heterogeneous decomposition. IEEE Trans. on Sustain. Ener., 2019.
10. L. Deng, Z. Li, H. Sun, Q. Guo, et al. Generalized locational marginal pricing in a heat-and-electricity-integrated market. IEEE Trans. on Smart Grid, 2019,10(6): 6414-6425.
11. L. Deng, X. Zhang, H. Sun. Real-time Autonomous Trading in the Electricity-and-Heat Distribution Market Based on Blockchain. 2019 IEEE Power Energy Society General Meeting, 2019: 1-5.
12. GM. Kopanos, MC. Georgiadis, EN. Pistikopoulos. Energy production planning of a network of micro combined heat and power generators. Appl. Energy, 2013, 102: 1522-1534.
13. A. Zidan, HA. Gabbar, A. Eldessouky. Optimal planning of combined heat and power systems within microgrids. Energy, 2015, 93: 235-244.
14. A. Gupte. Mixed integer bilinear programming with applications to the pooling problem. Diss. Georgia Institute of Technology, 2012.
15. Z. Li, W. Wu, M. Shahidehpour, J. Wang, B. Zhang. Combined heat and power dispatch considering pipeline energy storage of district heating network. IEEE Trans. on Sustain. Ener., 2015, 7(1): 12-22.
16. Y. Chen, Q. Guo, H. Sun, Z. Li, Z. Pan, W. Wu. A water mass method and its application to integrated heat and electricity dispatch considering thermal inertias. Energy, 2019, 181: 840-852.
17. S. Huang, W. Tang, Q. Wu, C. Li. Network constrained economic dispatch of integrated heat and electricity systems through mixed integer conic programming. Energy, 2019, 179: 464-474.
18. Y. Jiang, C. Wan, A. Botterud, Y. Song, M. Shahidehpour. Convex Relaxation of Combined Heat and Power Dispatch. arXiv preprint arXiv:2006.04028 (2020).
19. PM. Castro. Tightening piecewise McCormick relaxations for bilinear problems. Computers & Chemical Engineering, 2015, 72: 300-311.
20. L. Gurobi Optimization. Gurobi Optimizer Reference Manual. 2020, http://www.gurobi.com.
21. S. Chen, AJ. Conejo, R.Sioshansi, Z. Wei. Unit commitment with an enhanced natural gas-flow model. IEEE Trans. on Power Syst., 2019, 34(5): 3729-3738.
22. Test Data for the Day-ahead Planning Problem of Integrated Heat and Electricity Systems, 2020, https://docs.google.com/spreadsheets/d/1px21mm7BfATtTo11SStwiEEBSSpwy4JAP6Sn7O6TX_U/edit?usp=sharing.